\documentclass[reqno]{article}
\usepackage{amsfonts, amsmath}
\begin{document}
\setlength{\textwidth}{16cm}
\renewcommand{\thesection}{\arabic{section}.}
\renewcommand{\theequation}{\arabic{section}.\arabic{equation}}
\newcommand{\be}{\begin{eqnarray}}
\newcommand{\en}{\end{eqnarray}}
\newcommand{\no}{\nonumber}
\newcommand{\la}{\lambda}
\newcommand{\ep}{\epsilon}
\newcommand{\lan}{\langle}
\newcommand{\ra}{\rangle}
\newcommand{\de}{\delta}
\newcommand{\ov}{\overline}
\newcommand{\bet}{\beta}
\newcommand{\al}{\alpha}
\newcommand{\fr}{\frac}
\newcommand{\we}{\wedge}
\newcommand{\D}{\Delta}
\newcommand{\R}{\mathbb{R}}
\newcommand{\ri}{\rightarrow}
\newcommand{\om}{\Omega}
\newcommand{\na}{\nabla}
\newcommand{\vs}{\vskip0.3cm}
\newcommand{\pa}{\partial}
\newcommand{\va}{\varsigma}
\pagestyle{myheadings}
\newtheorem{theorem}{Theorem}[section]
\newtheorem{lemma}[theorem]{Lemma}
\newtheorem{condition}[theorem]{Conditions}
\newtheorem{corollary}[theorem]{Corollary}
\newtheorem{proposition}[theorem]{Proposition}
\newtheorem{remark}[theorem]{Remark}
\newtheorem{example}[theorem]{Example}
\newtheorem{conjecture}[theorem]{Conjecture}
\newtheorem{hypotheses}[theorem]{Hypotheses}
\newtheorem{defn}[theorem]{Definition}
\newtheorem{maintheorem}{Main Theorem}[section]
\renewcommand{\thefootnote}{}

\title{On a Conjecture of Ashbaugh and Benguria about Lower   Eigenvalues of the Neumann Laplacian}
\footnotetext{ 
\\
\ \ \ \ \ \hspace*{4ex} Key Words: Isoperimetric Inequality,
Eigenvalues, Free Membrane Problem, Szeg\"o-Weinberger Inequality,
Ashbaugh-Benguria Conjecture.}
\author{
Qiaoling Wang$^{a}$, Changyu Xia$^{a}$}

\date{}

\maketitle ~~~\\[-15mm]

\begin{center}
{\footnotesize  $a$. Departamento de Matem\'atica, Universidade de
Brasilia, 70910-900-Brasilia-DF, Brazil Email: wang@mat.unb.br(Q.
Wang), xia@mat.unb.br(C. Xia).  }
\end{center}

\begin{abstract}
In this paper, we prove an isoperimetric inequality for lower order
eigenvalues of the free membrane problem on bounded domains in a
Euclidean space or a hyperbolic space which strengthens the
well-known Szeg\"o-Weinberger inequality and supports  a celebrated
conjecture of Ashbaugh-Benguria.
 \end{abstract}

\markright{\sl\hfill  Q. Wang, C. Xia \hfill}

\section{Introduction}
\renewcommand{\thesection}{\arabic{section}}
\renewcommand{\theequation}{\thesection.\arabic{equation}}
\setcounter{equation}{0}

Let $(M, g)$ be complete Riemannian manifold of dimension $n$,
$n\geq 2$. We denote by $\D $ the Laplace operator on $M$. For
bounded domain  $\om$ with smooth boundary  in $M$ we consider the
free membrane problem
 \be\label{1.1}\left\{\begin{array}{l}
\Delta f = -\mu f \ \ \ {\rm in  \ }\ \ \om,\\
\fr{\pa f}{\pa\nu} =  0 \ \ \ \ \ \ \ \ {\rm on  \ }\ \pa\om.
\end{array}\right.
\en Here  $\fr{\pa }{\pa\nu}$  denotes the outward unit normal
derivative on $\pa\om$. It is well known that the problem
(\ref{1.1}) has  discrete spectrum consisting in a sequence \be \no
\mu_0=0<\mu_1\leq \mu_2\leq\cdots \rightarrow +\infty.\en In the two
dimensional case,  G. Szego \cite{s} proved via conformal mapping
techniques that if $\om\subset \R^2 $ is simply connected, then
\be\label{in1} \mu_1(\om)A(\om)\leq \left.(\mu_1 A)\right|_{disk}=
\pi p_{1,1}^2 \en  where $A$ denotes the area. Later, using more
general methods,
 Weinberger \cite{w} showed that (\ref{in1}) and its
$n$-dimensional analogue, \be \label{int2} \mu_1(\om)\leq
\left(\fr{\omega_n}{|\om|}\right)^{2/n} p_{n/2,1}^2,\en
 hold for
arbitrary domains in $\R^2$ and $\R^n$, respectively. Here $J_v$ is
the Bessel function of the first kind of order $v$, $p_{v,k}$ is the
$k$th positive zero of the derivative of $x^{l-v}J_v(x)$ and $|\om|$
denotes the volume of $\om$. Szeg\"o and Weinberger also noticed
that Szeg\"o's proof of (\ref{in1}) for simply connected domains in
$\R^2$ extends to prove the bound \be\label{in3}\fr 1{\mu_1}+\fr
1{\mu_2}\geq \fr{2A}{\pi p_{1,1}^2}\en for such domains. The bounds
of Szeg\"o and Weinberger are isoperimetric with equality if and
only if $\om$ is a disk ($n$-dimensional ball in the case of
Weinberger's result (\ref{int2})). A quantitative improvement of
(\ref{in1}) was made by  Brasco and Pratelli in \cite{bp} who showed
that for any bounded domain with smooth boundary $\om\subset\R^n$ we
have \be \omega_{n}^{2/n}p_{n/2,1}^2 -\mu_1(\om)|\om|^{2/n}\geq
c(n)\mathcal{A}(\om)^2.\en Here, $c(n)$ is positive constant
depending only on $n$ and  ${\cal A}(\om)$ is the so called {\it
Fraenkel asymmetry}, defined by $$ {\cal
A}(\om)=\inf\left\{\fr{|\om\Delta B|}{|\om|}: \ B\ {\rm ball\ in \
\mathbb{R}^n\ such\ that\ } |B|=|\om|\right\}. $$ Nadirashvilli
obtained in \cite{n} a  quantitative improvement of (\ref{in3})
which states that
 there exists a constant $C>0$ such that for every $\om\subset\R^2$
 smooth simply connected bounded open set it holds
 \be
 \fr 1{|\om|}\left(\fr 1{\mu_1(\om)}+\fr 1{\mu_2(\om)}\right)
 -\fr 1{|B|}\left(\fr 1{\mu_1(B)}+\fr 1{\mu_2(B)}\right)\geq \fr
 1C{\cal A}(\om)^2,
 \en
 where $B$ is any disk in $\R^2$. On the other hand,  Ashbaugh and Benguria
\cite{ab1} showed that \be \label{int3} \fr 1{\mu_1(\om)}+\cdots
+\cdots \fr 1{\mu_n(\om)} \geq \fr
n{n+2}\left(\fr{|\om|}{\omega_n}\right)^{2/n} \en holds for any
$\om\subset\R^n$. Some generalizations to (1.7) haven been obtained
e.g., in \cite{hx}, \cite{x}.

In\cite{ab1}, Ashbaugh and Benguria also proposed the following
important \vskip0.3cm {\bf Conjecture I }(\cite{ab1}). {\it For any
 bounded domain $\om$ with smooth boundary in $\R^n$, we have \be
\label{int3} \fr 1{\mu_1(\om)}+\fr 1{\mu_2(\om)}+\cdots \fr
1{\mu_n(\om)} \geq \fr{n\left(|\om|/\omega_n\right)^{2/n}}{
p_{n/2,1}^2} \en with equality holding if and only if $\om$ is a
ball in $\R^n$. }

\vskip0.3cm Ashbaugh \cite{a} and Henrot \cite{h} mentioned this
conjecture again.

A more general conjecture might be true. That is, for any bounded
domain domain $\om$ with smooth boundary in $\R^n$, it would
hold
\be\no \label{int2} \mu_n(\om)\leq
\left(\fr{\omega_n}{|\om|}\right)^{2/n} p_{n/2,1}^2,\en with
equality holding if and only if $\om$ is a ball in $\R^n$.

The Szeg\"o-Weinberger inequality (\ref{int2}) has been generalized
to bounded domains in a hyperbolic space  by Ashbaugh-Benguria
\cite{ab2} and Xu \cite{x} independently.  In his book, Chavel
\cite{c} mentioned that one can use Weinberger's method to prove
this result. In \cite{ab2}, Ashbaugh-Benguria also proved the
Szeg\"o-Weinberger inequality for bounded domains in a hemisphere.
One can also consider similar estimates for lower order eigenvalues
of the Neumann Laplacian for bounded domains in a hyperbolic space
or a hemisphere. \vskip0.3cm {\bf Conjecture II. } {\it Let $M$  be
an $n$-dimensional complete simply connected Riemannian manifold of
constant sectional curvature $\kappa\in\{-1, 1\}$  and $\om$ be a
bounded domain in $M$ which is contained in a hemisphere in the case
that $\kappa=1$. Let $B_{\om}$ be a geodesic ball in $M$ such that
$|\om|=|B_{\om}|$ and denote by $\mu_1(B_{\om})$ the first nonzero
eigenvalue of the Neumann Laplacian of $B_{\om}$. Then  the first
$n$ non-zero eigenvalues of the Neumann Laplacian of $\om$ satisfy
\be \label{conj2.1} \fr 1{\mu_1(\om)}+\fr 1{\mu_2(\om)}+\cdots \fr
1{\mu_n(\om)} \geq \fr{n}{\mu_1(B_{\om})} \en with equality holding
if and only if $\om$ is isometric to $B_{\om}$.} \vskip0.3cm

In this paper, we prove an isoperimetric inequality for the sums of
the reciprocals of the first $(n-1)$ non-zero eigenvalues of the
Neumann Laplacian on bounded domains in $\R^n$ or a hyperbolic space
which supports the above conjectures.

\begin{theorem}\label{th1} Let $\om$ be a bounded domain  with smooth boundary in
$\R^n$. Then \be \label{int3} \fr 1{\mu_1(\om)}+\cdots +\fr
1{\mu_{n-1}(\om)} \geq \fr{(n-1)\left(|\om|/\omega_n\right)^{2/n}}{
p_{n/2,1}^2} \en with equality holding if and only if $\om$ is a
ball in $\R^n$.
\end{theorem}

\begin{theorem}\label{th2}
Let $\mathbb{H}^n$  be an $n$-dimensional hyperbolic space of
curvature $-1$  and $\om$ be a bounded domain in $\mathbb{H}^n$. Let
$B_{\om}$ be a geodesic ball in $\mathbb{H}^n$ such that
$|\om|=|B_{\om}|$. Then we have \be \label{th2.1} \fr
1{\mu_1(\om)}+\cdots +\fr 1{\mu_{n-1}(\om)} \geq
\fr{n-1}{\mu_1(B_{\om})} \en with equality holding if and only if
$\om$ is isometric to $B_{\om}$.
\end{theorem}

\section{A proof of Theorem 1.1.}
\setcounter{equation}{0} In this section, we shall  prove the
following result which implies Theorem \ref{th1}.
\begin{theorem}\label{th2.1} Let $\om$ be a bounded domain  with smooth boundary in
$\R^n$. There exists a positive constant $d(n)$ depending only on
$n$ such that  the first $(n-1)$ nonzero Neumann eigenvalues of the
Laplacian of $\om$ satisfy the inequality \be
\omega_{n}^{2/n}p_{n/2,1}^2 -\fr{(n-1)|\om|^{2/n}}{\fr
1{\mu_1}+\cdots + \fr 1{\mu_{n-1}}}\geq d(n)\mathcal{A}(\om)^2,\en
with equality holding if and only if $\om$ is an $n$-ball.
\end{theorem}

{\bf Remark.} One can easily see that (2.1) strengthens (1.5).
\vskip0.2cm
 Before proving Theorem \ref{th2.1}, we recall some known
facts we need (Cf. \cite{c},\cite{h},\cite{sy}).

Let $\{u_j\}_{j=0}^{\infty}$ be an orthonormal set of
eigenfunctions of the problem (\ref{1.1}), that is,
\be\label{pth1.6} \left\{\begin{array}{l} \Delta u_i= -\mu_i u_i \ \
\ {\rm in} \ \ \
\om,\\
\left.\fr{\pa u_i}{\pa\nu}\right|_{\pa\om}=0,\\
\int_{\om} u_i u_j dv_g=\delta_{ij}.
\end{array}\right.,  \en where
$dv_g$ denotes the volume element of the metric $g$. For each
$i=1,2,\cdots,$ the variational characterization of $\mu_i(\om)$ is
given by
 \be\label{pth0.5}
\mu_i(\om)=\inf_{u\in
H^1(\om)\setminus\{0\}}\left\{\fr{\int_{\om}|\na u|^2
dv_g}{\int_{\om} u^2 dv_g}: \int_{\om} uu_j dv_g=0, j=0,\cdots,
i-1\right\}.
 \en
Let $B_r$ be a ball of radius $r$ centered at the origin in $\R^n$.
It is known that $\mu_1(B_r)$ has multiplicity $n$, that is,
$\mu_1(B_r)=\cdots =\mu_{n}(B_r)$. This value can be explicitly
computed together with its corresponding eigenfunctions.
 A basis for the eigenspace corresponding to $\mu_1(B_r)$
consists of \be\label{pth1.2} \xi_i(x)=|x|^{1-\fr n2}
J_{n/2}\left(\fr{p_{n/2, 1}|x|}r\right)\fr{x_i}{|x|}, \ \
i=1,\cdots, n. \en The radial part of $\xi_i$ \be\label{pth1.3}
g(|x|)=|x|^{1-\fr n2} J_{n/2}\left(\fr{p_{n/2, 1}|x|}r\right), \en
satisfies the differential equation of Bessel type
\be\label{pth1.4}\left\{\begin{array}{l}
g^{\prime\prime}(t)+\fr{n-1}t
g^{\prime}(t)+\left(\mu_1(B_r)-\fr{n-1}{t^2}\right)g(t)=0,\\
g(0)=0, \ \ g^{\prime}(r)=0.
\end{array}\right.
\en We can compute \be \label{pth1.52} \mu_1(B_r)&=&
\fr{\int_{B_r}\left( g^{\prime}(|x|)^2
+(n-1)\fr{g(|x|)^2}{|x|^2}\right)dx}{\int_{B_r}g(|x|)^2 dx}
\\ \no
&=& \left(\fr{p_{n/2, 1}}r\right)^2. \en {\it Proof of Theorem
\ref{th2.1}.}  Let \be\label{pth1.51}
r=\left(\fr{|\om|}{\omega_n}\right)^{1/n}\en
 and
define $G: [0, +\infty)\ri\R$ by \be\label{pth1.6}
G(t)=\left\{\begin{array}{l}
g(t), \ t\leq r,\\
g(r), \ t > r.
\end{array}
\right. \en   We need to choose suitable trial functions $\phi_i$
for each of the eigenfunctions $u_i$ and insure that these are
orthogonal to the preceding eigenfunctions $u_0,\cdots, u_{i-1}$.
For the $n$ trial functions $\phi_1, \phi_2, \cdots, \phi_n,$ we
choose:
 \be \phi_i= G(|x|)\fr{x_i}{|x|}, \ \ {\rm for}\
\ i=1,\cdots, n, \en but before we can use these we need to make
adjustments  so that \be \phi_i \perp{\rm span}\{u_0,\cdots,
u_{i-1}\} \en
 in $L^2(\om)$. In order to do this, let us fix an orthonormal basis
 $\{e_i\}_{i=1}^n$ of $\R^n$. From the well-know arguments of
 Weinberger in \cite{w}
by using the Brouwer fixed point theorem,
 we know that it is always possible to choose the origin of
$\R^n$ so that \be \label{pth1.7} \int_{\om}\langle x,
e_i\rangle\fr{G(|x|)}{|x|} dx =0,\ \ i=1,\cdots, n,\en that is, $
\langle x, e_i\rangle\fr{G(|x|)}{|x|}\perp u_0$ (which is actually
just the constant function $1/\sqrt{|\om|}$). Here  $dx$ and
$\langle , \rangle$ denote the standard Lebesgue measure and the
inner product of $\R^n$, respectively.   Now we show that there
exists a new orthonormal basis $\{e_i^{\prime}\}_{i=1}^n$ of $\R^n$
such that \be\label{pth1.8} \langle x,
e_i^{\prime}\rangle\fr{G(|x|)}{|x|}\perp u_j, \en for $j=1,\cdots,
i-1$ and $i=2,\cdots, n$. To see this, we define an $n \times n$
matrix $Q=\left(q_{ij}\right)$ by  \be q_{ij}=\int_{\om} \langle x,
e_i\rangle\fr{G(|x|)}{|x|} u_j(x) dx, \ i,j=1,2,\cdots,n.\en Using
the orthogonalization of Gram and Schmidt (QR-factorization
theorem), we know that there exist an upper triangle matrix
$T=(T_{ij})$ and an orthogonal matrix $U=(a_{ij})$ such that $T=UQ$,
i.e.,
\begin{eqnarray*}
T_{ij}=\sum_{k=1}^n a_{ik}q_{kj}=\int_{\om} \sum_{k=1}^n
a_{ik}\langle x, e_k\rangle\fr{G(|x|)}{|x|} u_j(x) dx =0,\ \ 1\leq
j<i\leq n.
\end{eqnarray*}
Letting $e_i^{\prime}=\sum_{k=1}^n  a_{ik}e_k, \ i=1,...,n$;  we
arrive at (\ref{pth1.8}). Let us denote by $x_1, x_2,\cdots, x_n$
the coordinate functions with respect to the base
$\{e_i^{\prime}\}_{i=1}^n$, that is, $x_i=\langle x,
e_i^{\prime}\rangle, \ x\in\R^n$. From (\ref{pth1.7}) and
(\ref{pth1.8}), we have \be\label{pth1.9}\int_{\om}\phi_i u_j dx=
\int_{\om} G(|x|)\fr{x_i}{|x|} u_j(x) dx=0, \ i=1,\cdots, n, \
j=0,\cdots, i-1. \en It then follows from the variational
characterization (\ref{pth0.5}) that
\begin{eqnarray}\label{pth1.11}
\mu_{i}\int_{\om} \phi_i^2dx \leq\int_{\om} |\na\phi_i|^2 dx,\
i=1,\cdots,n. \en Substituting  \be\label{p} |\na \phi_i|^2&=&
G^{\prime}(|x|)^2\fr{x_i^2}{|x|^2}+\fr{G(|x|)^2}{|x|^2}\left(1-\fr{x_i^2}{|x|^2}\right)\\
\no &=&
\fr{G(|x|)^2}{|x|^2}+\left(G^{\prime}(|x|)^2-\fr{G(|x|)^2}{|x|^2}\right)\fr{x_i^2}{|x|^2}
\en into (\ref{pth1.11}) and dividing by $\mu_i$, one gets for
$i=1,\cdots,n$ that \be\label{pth1.12} \int_{\om}\phi_i^2 dx\leq\fr
1{\mu_i}\int_{\om} \fr{G(|x|)^2}{|x|^2}dx + \fr 1{\mu_i}\int_{\om}
\left(G^{\prime}(|x|)^2-\fr{G(|x|)^2}{|x|^2}\right)\fr{x_i^2}{|x|^2}dx.
\en Summing over $i$, we get \be\label{d7}\int_{\om} G(|x|)^2 dx
 &\leq&  \sum_{i=1}^n \fr
1{\mu_i}\int_{\om}\fr{G(|x|)^2}{|x|^2}dx \\ \no & & +\sum_{i=1}^n
\fr 1{\mu_i}\int_{\om}
\left(G^{\prime}(|x|)^2-\fr{G(|x|)^2}{|x|^2}\right)\fr{x_i^2}{|x|^2}dx.
\en
 Since
\be \sum_{i=1}^n\fr 1{\mu_i}\fr{x_i^2}{|x|^2}&=& \sum_{i=1}^{n-1}\fr
1{\mu_i}\fr{x_i^2}{|x|^2}+ \fr 1{\mu_n} \fr{x_n^2}{|x|^2}\\ \no &=&
\sum_{i=1}^{n-1}\fr 1{\mu_i} \fr{x_i^2}{|x|^2}+ \fr
1{\mu_n}\left(1-\sum_{i=1}^{n-1} \fr{x_i ^2}{|x|^2}\right), \en we
have \be & &\sum_{i=1}^n \fr 1{\mu_i}\int_{\om}
\left(G^{\prime}(|x|)^2-\fr{G(|x|)^2}{|x|^2}\right)\fr{x_i^2}{|x|^2}dx
\\ \no &=& \sum_{i=1}^{n-1} \fr 1{\mu_i}\int_{\om}
\left(G^{\prime}(|x|)^2-\fr{G(|x|)^2}{|x|^2}\right)\fr{x_i^2}{|x|^2}dx
\\ \no & & + \fr 1{\mu_n}\int_{\om}
\left(G^{\prime}(|x|)^2-\fr{G(|x|)^2}{|x|^2}\right)dx \\ \no & & -
\fr 1{\mu_n}\int_{\om}
\left(G^{\prime}(|x|)^2-\fr{G(|x|)^2}{|x|^2}\right)\sum_{i=1}^{n-1}\fr{x_i^2}{|x|^2}
dx
\\ \no&=&
\sum_{i=1}^{n-1}\int_{\om}\left(\fr 1{\mu_i}-\fr
1{\mu_n}\right)\left(G^{\prime}(|x|)^2-\fr{G(|x|)^2}{|x|^2}\right)\fr{x_i^2}{|x|^2}dx
\\ \no & & + \fr 1{\mu_n}\int_{\om}
\left(G^{\prime}(|x|)^2-\fr{G(|x|)^2}{|x|^2}\right) dx.\en

\begin{lemma}\label{le1} We have $g^{\prime}|_{[0, r)}>0$, $g|_{(0, r]}>0$ and $g^{\prime}(t)-\fr{g(t)} t\leq 0, \ \forall t\in (0, r].$
\end{lemma}
{\it Proof of Lemma \ref{le1}.}
 The Bessel function of the first
kind $J_v(t)$ is given by \be
J_v(t)=\sum_{k=0}^{+\infty}\fr{(-1)^k\left(\fr
t2\right)^{2k+v}}{k!\Gamma(k+v+1)}, \en which, combining with
(\ref{pth1.3}), gives \be\label{le1.2} g(t) =\left(\fr{p_{n/2,
1}}{2r}\right)^{\fr
n2}t\sum_{k=0}^{+\infty}\fr{(-1)^k\left(\fr{p_{n/2,
1}}{2r}t\right)^{2k}}{k!\Gamma(k+\fr n2 +1)}. \en Thus, $g(0)=0,
g^{\prime}(0)>0$.  Since $r$ is the first positive zero of
$g^{\prime}$, we have $g|_{(0, r]}>0$ and $g^{\prime}|_{[0, r)}>0.$
Observe that \be\label{le1.21}
  \ \lim_{t\ri
0}\left(g^{\prime}(t)-\fr{g(t)}t\right)=0, \ \
g^{\prime}(r)-\fr{g(r)}r<0. \en Let us assume by contradiction that
there exists a $t_0\in (0, r)$ such that \be
g^{\prime}(t_0)-\fr{g(t_0)}{t_0}>0.\en In this case, we know from
(\ref{le1.21})  that the function $g^{\prime}(t)-\fr{g(t)}t$ attains
its maximum at some  $t_1\in (0, r)$ and so we have \be
\label{le1.4}
g^{\prime\prime}(t_1)-\fr{t_1g^{\prime}(t_1)-g(t_1)}{t_1^2}=0. \en
From (2.6), we have
 \be\label{le1.5}
g^{\prime\prime}(t_1)+\fr{n-1}{t_1}
g^{\prime}(t_1)+\left(\mu_1(B_r)-\fr{n-1}{t_1^2}\right)g(t_1)=0. \en
Eliminating $g^{\prime\prime}(t_1)$ from (\ref{le1.4}) and
(\ref{le1.5}), we get \be \fr
n{t_1}\left(g^{\prime}(t_1)-\fr{g(t_1)}{t_1}\right) =
-\mu_1(B_r)g(t_1)<0. \en This is a contradiction and completes the
proof of Lemma 2.1.

From Lemma 2.1 and the definition of $G$, we  know that \be
G^{\prime}(|x|)^2-\fr{G(|x|)^2}{|x|^2}\leq 0 \ \ \ {\rm on\ \ }\om.
\en Hence \be
 \sum_{i=1}^{n-1}\int_{\om}\left(\fr
1{\mu_i}-\fr
1{\mu_n}\right)\left(G^{\prime}(|x|)^2-\fr{G(|x|)^2}{|x|^2}\right)\fr{x_i^2}{|x|^2}dx\leq
0. \en Combining (2.19), (2.21) and (2.30), one gets
\be\label{d7}\int_{\om} G(|x|)^2 dx
 &\leq&  \fr 1{\mu_n}\int_{\om}
\left(G^{\prime}(|x|)^2-\fr{G(|x|)^2}{|x|^2}\right)dx\\
\no & & +\sum_{i=1}^{n} \fr 1{\mu_i}\int_{\om}\fr{G(|x|)^2}{|x|^2}dx
\\ \no &=& \fr 1{\mu_n}\int_{\om}
G^{\prime}(|x|)^2 +\sum_{i=1}^{n-1} \fr
1{\mu_i}\int_{\om}\fr{G(|x|)^2}{|x|^2}dx \\ \no &\leq& \fr 1{n-1}
\sum_{i=1}^{n-1}\fr 1{\mu_i}\int_{\om}\left(G^{\prime}(|x|)^2
+(n-1)\fr{G(|x|)^2}{|x|^2}\right)dx, \en that is, \be \fr
{n-1}{\sum_{i=1}^{n-1}\fr 1{\mu_i} }\int_{\om} G(|x|)^2 dx \leq
\int_{\om}\left(G^{\prime}(|x|)^2
+(n-1)\fr{G(|x|)^2}{|x|^2}\right)dx. \en Using the fact that $G(t)$
is increasing, one gets \be \label{th1.9} \int_{\om} G(|x|)^2dx &=&
\int_{\om\cap B_r} G(|x|)^2dx+\int_{\om\setminus B_r} G(|x|)^2dx\\
\no &\geq &
 \int_{\om\cap B_r}
G(|x|)^2dx+ g(r)^2|\om\setminus B_r|
\\ \no &=&
 \int_{\om\cap B_r}
g(|x|)^2dx+ g(r)^2|B_r\setminus \om|
\\ \no &\geq &
 \int_{\om\cap B_r}
g(|x|)^2+\int_{B_r\setminus \om} g(|x|)^2dx
\\ \no &=& \int_{B_r}g(|x|)^2dx,
\en which, combining with (2.32), gives \be \fr
{n-1}{\sum_{i=1}^{n-1}\fr 1{\mu_i} }\int_{B_r}g(|x|)^2dx \leq
\int_{\om}\left(G^{\prime}(|x|)^2
+(n-1)\fr{G(|x|)^2}{|x|^2}\right)dx. \en We know from (2.7) that \be
\left(\fr{p_{n/2,1}}{r}\right)^2\int_{B_r}g(|x|)^2dx &=&
\int_{B_r}\left(g^{\prime}(|x|)^2
+(n-1)\fr{g(|x|)^2}{|x|^2}\right)dx \\ \no &=&
\int_{B_r}\left(G^{\prime}(|x|)^2
+(n-1)\fr{G(|x|)^2}{|x|^2}\right)dx.\en Consequently, we have \be &
& \left(\left(\fr{p_{n/2,1}}{r}\right)^2-\fr
{n-1}{\sum_{i=1}^{n-1}\fr 1{\mu_i} } \right)\int_{B_r}g(|x|)^2dx \\
\no &\geq& \int_{B_r} \left(G^{\prime}(|x|)^2
+(n-1)\fr{G(|x|)^2}{|x|^2}\right)dx-\int_{\om}\left(G^{\prime}(|x|)^2
+(n-1)\fr{G(|x|)^2}{|x|^2}\right)dx. \en We have \be \no\fr
d{dt}\left[G^{\prime}(t)^2+(n-1)\fr{G(t)^2}{t^2}\right]=2G^{\prime}(t)G^{\prime\prime}(t)+2(n-1)(tG(t)G^{\prime}(t)-G(t)^2)/t^3.
\en For $t>r$ this is negative since $G$ is constant there. For
$t\leq r$ we use the differential equation (2.6) to obtain
 \be\no
\fr d{dt}\left[G^{\prime}(t)^2+(n-1)\fr{G(t)^2}{t^2}\right]=-2
\mu_1(B_r)GG^{\prime}-(n-1)(tG^{\prime}-G)^2/t^3<0.\en Thus the
function $ G^{\prime}(t)^2+(n-1)\fr{G(t)^2}{t^2}$ is decreasing for
$t>0$.

{\lemma \label{le2.2} (\cite{bp}) Let $f: \mathbb{R}_+\ri
\mathbb{R}_+$ be a decreasing function. Then we have \be
\label{le1.0} \int_{B_r} f(|x|)dx -\int_{\om} f(|x|)dx\geq n\omega_n
\int_{\rho_1}^{\rho_2} |f(t)-f(r)|t^{n-1} dt. \en Here
 \be
\label{le1.2} \rho_1=\left(\fr{|\om\cap B_r|}{\omega_n}\right)^{\fr
1n} \ \ \ {\rm and}\ \ \ \rho_2=\left(\fr{|\om|+|\om\setminus
B_r|}{\omega_n}\right)^{\fr 1n}.\en } Taking
$f(t)=G^{\prime}(t)^2+(n-1)\fr{G(t)^2}{t^2}$ in Lemma 2.3, we obtain
\be\no & & \int_{B_r} \left(G^{\prime}(|x|)^2
+(n-1)\fr{G(|x|)^2}{|x|^2}\right)dx-\int_{\om}\left(G^{\prime}(|x|)^2
+(n-1)\fr{G(|x|^2}{|x|^2}\right)dx\\ \no &\geq &
 n\omega_n
\int_{r}^{\rho_2} |f(t)-f(r)|t^{n-1} dt
\\  & =&
 n\omega_n
\int_{r}^{\rho_2} (f(r)-f(t))t^{n-1} dt \en Observe that \be
(f(r)-f(t))t^{n-1}= (n-1)g(r)^2\left(\fr 1{r^2}-\fr
1{t^2}\right)t^{n-1}, \ \ \ {\rm for\ } \rho_2\geq
 t\geq r.
 \en
Therefore, \be\no & &\int_{r}^{\rho_2} (f(r)-f(t))t^{n-1} dt\\
 & = & g(r)^2\cdot\left\{\begin{array}{l}
\fr{n-1}{nr^2}\left(\rho_2^n-r^n\right)-\fr{n-1}{n-2}\left(\rho_2^{n-2}-r^{n-2}\right), \ \ \ {\rm if\ \ } n>2,\\
\fr 1{2r^2}\left(\rho_2^2-r^2\right)-\ln\fr{\rho_2}r,  \ \ \ \ \ \ \
\ \ \ \ \ \ \ \ \ \ \ \ \ \ \ \ {\rm if\ \ } n=2.
\end{array}\right.
\en By using the definition of $\rho_2$ we have when $n>2$, \be
\rho_2^{n-2}-r^{n-2}&=& r^{n-2}\left[\left(1+\fr{|\om\setminus
B_r|}{|\om|}\right)^{\fr{n-2}n}-1\right]\\ \no &\leq&
r^{n-2}\left(\fr{n-2}n\fr{|\om\setminus
B_r|}{|\om|}-\fr{(n-2)2^{-\fr 2n-1}}{n^2}\left(\fr{|\om\setminus
B_r|}{|\om|}\right)^2\right), \en thanks to the elementary
inequality \be\no (1+t)^{\delta}\leq 1+\delta t +\fr{\delta(\delta
-1)}2\cdot 2^{\delta-2} t^2, \forall \ \delta\in (0, 1), \ \forall
t\in [0, 1], \en and when $n=2$, \be \ln\fr{\rho_2}r &=&\fr 12
\ln\left(1+\fr{|\om\setminus B_r|}{|\om|}\right)
\\ \no &\leq& \fr 12 \left(\fr{|\om\setminus
B_r|}{|\om|}-\fr 14 \left(\fr{|\om\setminus
B_r|}{|\om|}\right)^2\right), \en thanks to the elementary
inequality \be\no \ln(1+t)\leq t-\fr{t^2}4, \  \forall t\in [0, 1].
\en Since $|B_r|=|\om|$, we have $|\om\Delta B_r|=2|\om\setminus
B_r|$ and so $$ \fr{|\om\setminus B_r|}{|\om|}\geq \fr 12
\mathcal{A}(\om).$$ It then follows by  substituting (2.42) and
(2.43) into (2.41) that
 \be& & \int_{r}^{\rho_2} (f(r)-f(t))t^{n-1} dt\\
\no
 &=& g(r)^2\left(\fr{|\om\setminus B_r|}{|\om|}\right)^2\cdot\left\{\begin{array}{l}
r^{n-2}\cdot \fr{(n-1)2^{-\fr 2n -1}}{n^2}, \ \ \  {\rm if\ \ } n>2,\\
\fr 18, \ \ \ \ \ \ \ \ \ \ \ \ \ \ \ \ \ \ \ \ \ \ \ \ {\rm if\ \ }
n=2.
\end{array}\right.
\\
\no
 &\geq&\fr 14 g(r)^2\mathcal{A}(\om)^2\cdot\left\{\begin{array}{l}
r^{n-2}\cdot \fr{(n-1)2^{-\fr 2n -1}}{n^2}, \ \ \  {\rm if\ \ } n>2,\\
\fr 18, \ \ \ \ \ \ \ \ \ \ \ \ \ \ \ \ \ \ \ \ \ \ \ \ {\rm if\ \ }
n=2.
\end{array}\right.
\en Thus, concerning the right hand side of (2.36), one gets from
(2.39) and (2.44) that \be\no & & \int_{B_r} \left(G^{\prime}(|x|)^2
+(n-1)\fr{G(|x|)^2}{|x|^2}\right)dx-\int_{\om}\left(G^{\prime}(|x|)^2
+(n-1)\fr{G(|x|^2}{|x|^2}\right)dx\\ \no &\geq &
\fr{\omega_n}4g(r)^2\mathcal{A}(\om)^2\cdot \left\{\begin{array}{l}
r^{n-2}\cdot \fr{(n-1)2^{-\fr 2n
-1}}{n}, \ \ \  {\rm if\ \ } n>2,\\
\fr 14 , \ \ \ \ \ \ \ \ \ \ \ \ \ \ \ \ \ \ \ \ \ \ \ \  {\rm if\ \
} n=2,
\end{array}\right.
\\ \no &=& \fr{\omega_n}4 J_{n/2}(p_{n/2,1})^2\mathcal{A}(\om)^2\cdot
\left\{\begin{array}{l}\fr{(n-1)2^{-\fr 2n
-1}}{n}, \ \ \  {\rm if\ \ } n>2,\\
\fr 14, \ \ \ \ \ \ \ \ \ \ \ \ \ \ \ \ {\rm if\ \ } n=2,
\end{array}\right.
\\  &\equiv& \alpha(n)\mathcal{A}(\om)^2.\en
 Concerning the left hand side of (2.36), we have \be&
&\left(\left(\fr{p_{n/2,1}}r\right)^2-\fr
{n-1}{\sum_{i=1}^{n-1}\fr 1{\mu_i} } \right)\int_{B_r}g(|x|)^2dx\\
\no &=&\left(\left(\fr{p_{n/2,1}}r\right)^2-\fr
{n-1}{\sum_{i=1}^{n-1}\fr 1{\mu_i} } \right) r^2\int_{\{|y|\leq
1\}}|y|^{2-n} J_{\fr n2}(p_{n/2,1}|y|)^2 dy \\
\no &=&\left(p_{n/2,1}^2\omega_n^{2/n}-\fr
{(n-1)|\om|^{2/n}}{\sum_{i=1}^{n-1}\fr 1{\mu_i} } \right)\beta(n),
\en where \be\no \beta(n)=\omega_n^{-2/n}\int_{\{|y|\leq
1\}}|y|^{2-n} J_{\fr n2}(p_{n/2,1}|y|)^2 dy. \en Combining (2.36),
(2.45) and (2.46), we obtain \be p_{n/2,1}^2\omega_n^{2/n}-\fr
{(n-1)|\om|^{2/n}}{\sum_{i=1}^{n-1}\fr 1{\mu_i} }\geq
\alpha(n)\beta(n)^{-1}\mathcal{A}(\om)^2\equiv
d(n)\mathcal{A}(\om)^2. \en
 Moreover, we can see that
equality holds in (2.47) only when $\om$ is a ball. This completes
the proof of Theorem \ref{th2.1}.

\section{A Proof of Theorem \ref{th2}}
\setcounter{equation}{0}
 In this section, we shall prove Theorem
\ref{th2}. Firstly, we list some important facts we need. About each
point $p\in \mathbb{H}^n$ there exists a coordinate system
$(t,\xi)\in [0, +\infty)\times\mathbb{S}^{n-1}$  relative to which
the Riemannian metric reads as  \be ds^2=dt^2 +\sinh^2t d\sigma^2,
\en where $d\sigma^2$ is the canonical metric on the
$(n-1)$-dimensional unit sphere $\mathbb{S}^{n-1}$.
\begin{lemma}(Cf. \cite{c}, \cite{x}). \label{le2} Let  $B(p, r)$ be a
geodesic ball of radius $r$ with center $p$ in $\mathbb{H}^n$. Then
the eigenfunction corresponding to the first nonzero eigenvalue
$\mu_1(B(p, r))$ of the Neumann problem on $B(p, r)$ must be
\be\label{le3.1} h(t, \xi)=f(t)\omega(\xi), \ \xi\in
\mathbb{S}^{n-1}, \en where $\omega(\xi)$ is an eigenfunction
corresponding to the first nonzero eigenvalue of $\mathbb{S}^{n-1}$,
$f$ satisfies
\be\label{le2.2}\left\{\begin{array}{l}f^{\prime\prime}+(n-1)\coth
t+\left(\mu_1(B(p,r))-\fr{n-1}{\sinh^2 t}\right)f=0,
\\
f(0)=f^{\prime}(r)=0, \ f^{\prime}|_{[0, r)}\neq 0,
\end{array}
\right. \en and \be\label{le2.3}
\mu_1(B(p,r))=\fr{\int_{B(p,r)}\left(f^{\prime}(t)^2+(n-1)\fr{f(t)^2}{\sinh^2t}\right)dv}{\int_{B(p,r)}f(t)^2dv}.
\en

\end{lemma}
{\it Proof of Theorem \ref{th2}.} Assume that the radius of
$B_{\om}$ is $r$.  Let $f$ be as in Lemma 3.1. Noticing $f(t)\neq 0$
when $0<t\leq r$, we may assume that $f(t)>0$ for $0<t\leq r$ and so
$f$ is nondecreasing on $[0, r]$. Let $\{{\bf e}_i\}_{i=1}^n$ be an
orthonormal basis of $\mathbb{R}^n$ and set
$\omega_i(\xi)=\langle{\bf e}_i, \xi\rangle, \ \xi\in
\mathbb{S}^{n-1}\subset \mathbb{R}^n$. Define \be
F(t)=\left\{\begin{array}{l} f(t), \ \
t\leq r,\\
f(r), \ \ t>r.
\end{array}
\right. \en  Let us take  a point $p\in \mathbb{H}^n$ such that in
the above coordinate system at $p$ we have \be\label{pth2.1}
\int_{\om} F(t)\omega_i(\xi)dv=0, \ \ i=1,\cdots,n. \en Here, $dv$
is the volume element of $\mathbb{H}^n$. By using the same arguments
as in the proof of Theorem 2.1, we can  assume further that
\be\label{pth2.2} \int_{\om}F(t)\omega_i(\xi) u_jdv=0, \en for
$i=2,3,\cdots, n$ and $j=1,\cdots,i-1$. Here
$\{u_i\}_{i=0}^{+\infty}$ is a orthonormal set of eigenfunctions
corresponding to the eigenvalues $\{\mu_i(\om)\}_{i=0}^{+\infty}$.
Hence, we conclude from the Rayleigh-Ritz variational
characterization (\ref{pth0.5}) that \be & & \mu_i(\om)\int_{\om}
F(t)^2 \omega_i^2(\xi)dv\\ \no &\leq&
\int_{\om}|\nabla(F(t)\omega_i(\xi))|^2 dv\\ \no &=&
\int_{\om}\left(|F^{\prime}(t)|^2\omega_i^2(\xi)+F^2(t)|\tilde{\nabla}\omega_i(\xi)|^2\sinh^{-2}
t\right)dv, \ i=1,\cdots,n, \en where $\tilde{\nabla}$ denotes the
gradient
operator of $\mathbb{S}^{n-1}.$ Thus \be\label{pth2.3} & & \int_{\om} F(t)^2 \omega_i^2(\xi)dv\\
\no &\leq& \fr
1{\mu_i(\om)}\int_{\om}|F^{\prime}(t)|^2\omega_i^2(\xi)dv+ \fr
1{\mu_i(\om)}\int_{\om}F^2(t)|\tilde{\nabla}\omega_i(\xi)|^2\sinh^{-2}
t dv. \en Observing $F^{\prime}(t)=0, t\geq r$, one gets
\be\label{pth2.4}
\int_{\om}|F^{\prime}(t)|^2\omega_i^2(\xi)dv&=&\int_{\om\cap
B(p,r)}|F^{\prime}(t)|^2\omega_i^2(\xi)dv\\ \no &\leq& \int_{
B(p,r)}|F^{\prime}(t)|^2\omega_i^2(\xi)dv\\ \no &=& \int_0^r
\int_{\mathbb{S}^{n-1}}|F^{\prime}(t)|^2\omega_i^2(\xi)\sinh^{n-1} t
dA\ dt
\\ \no &=& \fr 1n \int_0^r
\int_{\mathbb{S}^{n-1}}|F^{\prime}(t)|^2\sinh^{n-1} t dA\ dt
\\ \no &=& \fr 1n\int_{B(p,r)}|F^{\prime}(t)|^2 dv,
\en where $dA$ denotes the area element of $\mathbb{S}^{n-1}$. Since
\be |\tilde{\nabla}\omega_i(\xi)|\leq 1, \ \ \sum_{i=1}^n
|\tilde{\nabla}\omega_i(\xi)|^2=n-1, \en we have \be\label{pth2.5} &
& \sum_{i=1}^n \fr
1{\mu_i(\om)}|\tilde{\nabla}\omega_i(\xi)|^2\\
\no &=&\sum_{i=1}^{n-1} \fr
1{\mu_i(\om)}|\tilde{\nabla}\omega_i(\xi)|^2 + \fr
1{\mu_n(\om)}\sum_{i=1}^{n-1}\left(1-|\tilde{\nabla}\omega_i(\xi)|^2\right)
\\
\no &\leq&\sum_{i=1}^{n-1} \fr
1{\mu_i(\om)}|\tilde{\nabla}\omega_i(\xi)|^2 + \sum_{i=1}^{n-1}\fr
1{\mu_i(\om)}\left(1-|\tilde{\nabla}\omega_i(\xi)|^2\right)\\
\no &=& \sum_{i=1}^{n-1}\fr 1{\mu_i(\om)}.\en Summing on $i$ from
$1$ to $n$ in (\ref{pth2.3}) and using (\ref{pth2.4}) and
(\ref{pth2.5}), we get \be\label{pth2.6} & & \int_{\om} F(t)^2dv\\
\no &\leq& \sum_{i=1}^n \fr
1{n\mu_i(\om)}\int_{B(p,r)}|F^{\prime}(t)|^2dv+\sum_{i=1}^{n-1}\fr
1{\mu_i(\om)}\int_{\om}F^2(t)\sinh^{-2}tdv. \en We need the
following lemma.
\begin{lemma}\label{le2}The function $h(t)=\fr{F(t)}{\sinh t}$  is decreasing.
\end{lemma}
{\it Proof of Lemma \ref{le2}.} Observe that \be\no \lim_{t\ri 0}
h(t)= f^{\prime}(0). \en
 Let us show that \be\label{ple2.1}
\gamma(t)\equiv f^{\prime}(t)-\coth t f(t)\leq 0, \ t\in (0, r]. \en
Since \be\lim_{t\ri 0} \gamma(t)=0, \ \gamma(r)= -\coth r f(r)<0,
\en if $\gamma(t_0)>0$ for some $t_0\in (0, r),$ then $\gamma$
attains its maximum at some $t_1\in (0, r)$ and so \be\label{ple2.2}
0=\gamma^{\prime}(t_1)=
f^{\prime\prime}(t_1)+\fr{f(t_1)}{\sinh^2t_1}-\coth t_1
f^{\prime}(t_1). \en We have from (\ref{le2.2}) that
\be\label{ple2.3} f^{\prime\prime}(t_1)+(n-1)\coth
t_1f^{\prime}(t_1)+ \mu_1(B(r))f(t_1)-\fr{n-1}{\sinh^2 t_1}f(t_1)=0.
\en Hence \be  f^{\prime}(t_1)-\fr{f(t_1)}{\cosh t_1\sinh t_1} =
-\fr{\mu_1(B(r))f(t_1)\sinh t_1}{n\cosh t_1}<0, \en which
contradicts to  \be f^{\prime}(t_1)-\coth t_1 f(t_1)>0. \en Thus
(\ref{ple2.1}) holds. Consequently $h^{\prime}(t)\leq 0, \ \forall
t\in (0, r]$ and $h$ is decreasing. The proof of Lemma 3.2 is
completed.

\vskip0.2cm Now we  go on the proof of Theorem \ref{th2}. Since $F$
is increasing and $\fr{F(t)}{\sinh t}$ is decreasing, we can use the
same arguments as in  the proof of (2.33)  to conclude that
\be\label{pth2.7} \int_{\om} F(t)^2 dv\geq \int_{B(p,r)}f(t)^2 dv\en
and \be\label{pth2.8} \int_{\om} \fr{F(t)^2}{\sinh^2t} dv\leq
\int_{B(p,r)} \fr{f(t)^2}{\sinh^2t} dv. \en Substituting
(\ref{pth2.7}) and (\ref{pth2.8}) into (\ref{pth2.6}), one gets
\be\label{pth2.9} \fr 1{n-1}\sum_{i=1}^{n-1}\fr
1{\mu_i(\om)}&\geq&\fr{\int_{B(p,r)}f(t)^2dv}{\int_{B(p,r)}\left(f^{\prime}(t)^2+(n-1)\fr{f(t)^2}{\sinh^2t}\right)dv}\\
\no &=&\fr 1{\mu_1(B(p,r))} \en and equality holds if and only if
$\om=B(p,r)$. This completes the proof of Theorem \ref{th2}.

\section*{Acknowledgments}

 Q. Wang was partially supported by CNPq,
Brazil (Grant No. 307089/2014-2). C. Xia was partially supported by
CNPq, Brazil (Grant No. 306146/2014-2).


\begin{thebibliography}{9999}

\bibitem{a} M. S. Ashbaugh, Open problems on eigenvalues of the
Laplacian, Analytic and Geometric Inequalities and Applications.
(1999), 13-28.




\bibitem{ab1}  M. S. Ashbaugh and R. D. Benguria, Universal bounds for the low eigenvalues of Neumann Laplacians in
$N$ dimensions. Siam J. Math. Anal. {\bf 24} (1993), 557-570.

\bibitem{ab2} M. S. Ashbaugh and R. D. Benguria, Sharp upper bound to the first
nonzero Neumann eigenvalue for bounded domains in spaces of constant
curvature, J. London Math. Soc. 52 (1995) 402-416.

\bibitem{bd} L. Brasco and G. De Philippis, Spectral inequalities in
quantitative form. Shape optimization and spectral theory, 201-281,
De Gruyter Open, Warsaw, 2017.


\bibitem{bp} L. Brasco, A. Pratelli, Sharp stability of some spectral inequalities, Geom. Funct. Anal., 22 (2012), 107-135.

\bibitem{c} I. Chavel,  Eigenvalues in Riemannian geometry (Academic, New York, 1984).

\bibitem{h} A. Henrot, Extremum problems for eigenvalues of elliptic
operators,  Birkh\"auser Verlag, Basel-Boston-Berlin, x + 202 pp.,
2006.  ISBN 978-3-76437705-2.

\bibitem{hx} G. N. Hile, Z. Xu, Inequalities for sums of the
reciprocals of eigenvalues, J. Math. Anal. Appl. {\bf 180} (1993)
412-430.

\bibitem{n} N. Nadirashvili, Conformal maps and isoperimetric inequalities for eigenvalues of the Neumann problem. Proceedings
 of the Ashkelon Workshop on Complex Function Theory (1996), 197-201, Israel Math. Conf. Proc. 11, Bar-Ilan Univ., Ramat Gan, 1997.


\bibitem{sy} R. Schoen R, S. T. Yau,   Lectures on Differential Geometry, Cambridge, 2004, MA: International Press.


\bibitem{s} G. Szeg\"o, Inequalities for certain eigenvalues of a membrane of
given area, J. Rational Mech. Anal. 3 (1954) 343-356.

\bibitem{w} H. F. Weinberger, An isoperimetric inequality for the
$n$-dimensional free membrane problem,  J. Rational Mech. Anal. 5
(1956) 633-636.

\bibitem{x} C. Xia, A universal bound for the low eigenvalues of Neumann Laplacians on
compact domains in a Hadamard manifold. Monatsh. Math. 128 (1999),
165-171.
\bibitem{x} Y. Xu, The first nonzero eigenvalue of neumann problem
on Riemannian manifolds, J. Geom. Anal. 5 (1995), 151-165.







\end{thebibliography}
\end{document}